\documentclass[12pt]{article}

\usepackage{graphicx,tikz,color,url}
\usepackage{amsmath,amssymb,latexsym}
\usepackage{pgfplots}
\usepackage{mathrsfs}
\usepackage{cite}
\usepackage{soul}
\usetikzlibrary{arrows}

\newtheorem{theorem}{Theorem}[section]

\newtheorem{lemma}[theorem]{Lemma}
\newtheorem{proposition}[theorem]{Proposition}
\newtheorem{corollary}[theorem]{Corollary}

\newtheorem{observation}[theorem]{Observation}

\newcommand{\proof}{\noindent{\bf Proof.\ }}
\newcommand{\qed}{\hfill $\square$ \bigskip}

\newcommand{\diam}{{\rm diam}}
\newcommand{\Spch}{\chi_{S}}
\newcommand{\smallqed}{{\tiny ($\Box$)}}
\newcommand{\cS}{{\cal S}}

\begin{document}

\title{S-packing chromatic critical graphs}

\author{
G\"{u}lnaz Boruzanl{\i} Ekinci $^{a,}$\thanks{Email: \texttt{gulnaz.boruzanli@ege.edu.tr}}
     \and Csilla Bujt\'{a}s $^{b,c,}$\thanks{Email: \texttt{csilla.bujtas@fmf.uni-lj.si}}
      \and Didem Gözüpek $^{d,}$\thanks{Email: \texttt{didem.gozupek@gtu.edu.tr}}
	\and Sandi Klav\v zar $^{b,c,e,}$\thanks{Email: \texttt{sandi.klavzar@fmf.uni-lj.si}}
}
\maketitle

\begin{center}
$^a$ Department of Mathematics, Faculty of Science, Ege University, İzmir, Türkiye\\
	\medskip
	$^b$ Faculty of Mathematics and Physics, University of Ljubljana, Slovenia\\
	\medskip
	
	$^c$ Institute of Mathematics, Physics and Mechanics, Ljubljana, Slovenia\\
	\medskip
    $^d$ Department of Computer Engineering, Gebze Technical University, Türkiye
	\medskip

    $^e$ Faculty of Natural Sciences and Mathematics, University of Maribor, Slovenia

\end{center}

\begin{abstract}
For a non-decreasing sequence of positive integers $S=(s_1,s_2,\ldots)$, the $S$-packing chromatic number of a graph $G$ is denoted by $\chi_S(G)$. In this paper, $\chi_S$-critical graphs are introduced as the graphs $G$ such that $\chi_S(H) < \chi_S(G)$ for each proper subgraph $H$ of $G$. Several families of $\chi_S$-critical graphs are constructed, and $2$- and $3$-colorable $\chi_S$-critical graphs are presented for all packing sequences $S$, while $4$-colorable $\chi_S$-critical graphs are found for most of $S$. Cycles which are $\chi_S$-critical are characterized under different conditions. It is proved that for any graph $G$ and any edge $e \in E(G)$, the inequality $\chi_S(G - e) \ge \chi_S(G)/2$ holds. Moreover, in several important cases, this bound can be improved to $\chi_S(G - e) \ge (\chi_S(G)+1)/2$. The sharpness of the bounds is also discussed. Along the way an earlier result on $\chi_S$-vertex-critical graphs is supplemented.
\end{abstract}


\noindent
{\bf AMS Subj.\ Class.\ (2020)}: 05C15, 05C12

\section{Introduction}

Let $S=(s_1,s_2,\ldots)$ be a non-decreasing sequence of positive integers and let $G = (V(G), E(G))$ be a graph. A mapping $c: V(G)\rightarrow[k] = \{1,\ldots, k\}$ is an {\em $S$-packing $k$-coloring} of $G$ if the equality $c(u)=c(v)=i$ for $u\neq v\in V(G)$ implies $d_G(u,v)>s_i$. The {\em $S$-packing chromatic number}  $\chi_S(G)$ of $G$ is the smallest integer $k$ such that $G$ admits an $S$-packing $k$-coloring~\cite{goddard-2012}.

In the special case when $S = (1,2,3,\ldots)$, the $S$-packing chromatic number is the standard packing chromatic number $\chi_\rho$, which was first explored under the name broadcast chromatic number~\cite{Goddard} and given the present name in~\cite{Bresar}. The 2020 review article~\cite{Bresar2} on packing colorings (including $S$-packing colorings) contains 68 references, but research continues, see~\cite{bidine-2023, ferme-2021, Greedy, gregor-2024, grochowski-2025}. The greatest emphasis in recent years has been on $S$-packing colorings, especially on subcubic graphs, see~\cite{bresar-2025, bresar-2025b, elzain-2025+, holub-2023, kostochka-2021, liu-2020, mortada-2024, mortada-2025, yang-2023}. 

It should be stressed that the concept of $S$-packing coloring is very general. As said, it contains the packing coloring as a particular instance. In addition, the special case $S = (k,k,k,\ldots)$, $k\ge 1$, is studied in the literature as {\em $k$-distance colorings}, the corresponding $(k,k,k,\ldots)$-packing chromatic number is denoted by $\chi_k$. Note that $\chi_1 = \chi$. Up to 2008, these investigations were surveyed in~\cite{kramer-2008}, while for some recent related papers see~\cite{hasanvard-2024, koley-2024, la-2025}. In the last years, however, the main focus was on $2$-distance colorings of planar graphs, cf.~\cite{aoki-2025, chen-2025, deniz-2025, la-2024, zhu-2023}.

Independently, and almost simultaneously, two different packing criticality concepts were introduced. In~\cite{klavzar-2019}, a graph $G$ was defined to be {\em $\chi_\rho$-vertex-critical} if $\chi_\rho(G-u) < \chi_\rho(G)$ for each $u\in V(G)$. Moreover, if $G$ is a $\chi_\rho$-vertex-critical graph with $\chi_\rho(G) = k$, then $G$ is called {\em $k$-$\chi_\rho$-vertex-critical}. On the other hand, according to~\cite{bresar-2022}, $G$ is {\em $\chi_\rho$-critical} if $\chi_\rho(H) < \chi_\rho(G)$ for each proper subgraph $H$ of $G$. If $G$ has no isolated vertices, this is equivalent to the requirement that $\chi_\rho(G-e) < \chi_\rho(G)$ holds for each $e\in E(G)$. If $G$ is a $\chi_\rho$-critical graph with $\chi_\rho(G) = k$, then $G$ is called {\em $k$-$\chi_\rho$-critical}. The paper~\cite{ferme-2022} further investigated $\chi_\rho$-vertex-critical graphs and provided a characterization of $4$-$\chi_\rho$-vertex-critical graphs.

In the same way as packing colorings extend to $S$-packing colorings, one can extend $\chi_\rho$-vertex-critical graphs and $\chi_\rho$-critical graphs to $\chi_S$-vertex-critical graphs and $\chi_S$-critical graphs. The first of these generalizations has been done in~\cite{holub-2020}, while in the follow-up paper~\cite{klavzar-2023} a characterization of 4-{$\chi_S$}-vertex-critical graphs for packing sequences with $s_1=1$ and $s_2\geq3$ is given. The second of these generalizations, that is, $\chi_S$-critical graphs, has not yet been studied, we fill this gap in this paper. We say that $G$ is {\em $\chi_S$-critical} if $\chi_S(H) < \chi_S(G)$ for each proper subgraph $H$ of $G$. If $G$ is a $\chi_S$-critical graph with $\chi_S(G) = k$, then $G$ is called {\em $k$-$\chi_S$-critical}.
Note that we do not consider the empty graph as a proper subgraph. Then, by our definition, the isolated vertex $K_1$ has no proper subgraph, and it is $1$-$\chi_S$-critical for every packing sequence $S$.

The paper is organized as follows. In the next section, we give some definitions, introduce useful notation, and present basic observations about $S$-packing critical graphs. In Section~\ref{sec:families}, some families of $\chi_S$-critical graphs are discussed. We determine $2$-$\chi_S$-critical and $3$-$\chi_S$-critical graphs for all packing sequences $S$, and determine $4$-$\chi_S$-critical graphs for most of $S$. For the case of $4$-$\chi_S$-critical graphs, we supplement an earlier result from the literature on $\chi_S$-vertex-critical graphs. In Section~\ref{sec:cycles} we characterize cycles which are $\chi_S$-critical under different conditions. In Section~\ref{sec:edge-removal} we consider the impact of edge removal on the $S$-packing chromatic number.  We prove that $\chi_S(G - e) \ge \chi_S(G)/2$ for any graph $G$ and any edge $e \in E(G)$, and that in several important cases the bound can be improved to $\chi_S(G - e) \ge (\chi_S(G)+1)/2$. For many $S$, infinitely many sharp examples are constructed. In the last section we identify several open problems for further research.

\section{Preliminaries}

Let $G=(V(G), E(G))$ be a graph. The open neighborhood $N_G(u)$ of $u$ in $G$ is the set of the neighbors of $u$. A support vertex of $G$ is a vertex adjacent to a leaf. The girth of $G$ is denoted by $g(G)$. If $G$ has no cycles, we set $g (G)= \infty$. As usual, $\alpha(G)$ is the independence number of $G$. The distance $d_G(u,v)$ between $u,v\in V(G)$ is the shortest-path distance. The diameter of $G$ is denoted by $\diam(G)$. A subgraph $H$ of $G$ is isometric, if for every two vertices $u,v\in V(H)$ we have $d_H(u,v) = d_G(u,v)$. The path, the cycle, and the complete graph of order $n$ are respectively denoted by $P_n$, $C_n$, and $K_n$, while the order of a graph $G$ will be denoted by $n(G)$. 

The set of all \textit{packing sequences} will be denoted by ${\cal S}$, that is, 
$${\cal S} = \{(s_1, s_2, \dots):\ 1\le s_1 \le s_2 \le \cdots\}.$$
For a given $S\in {\cal S}$ we will always assume that $S = (s_1, s_2, \dots)$. Unless stated otherwise, the packing sequences are considered to be infinite in this paper.

We will consider sets of packing sequences such that some of their first coordinates are fixed or bounded from below. Instead of introducing the notation in general, consider the following example. Assume we wish to consider the set of packing sequences $S=(s_1,s_2,s_3,\dots)$ with $s_1 = 1$, $s_2 = 3$, and $s_3 \ge 4$. Then we set 
$${\cal S}_{1,3,\overline{4}} = \{(s_1,s_2,s_3,\ldots):\ s_1 = 1, s_2 = 3, s_3 \ge 4\}\,.$$
The general notation should be clear from this example. For instance, ${\cal S}_{1,\overline{3},5}$ is the set of packing sequences with $s_1=1$, $3\le s_2\le 5$, and $s_3=5$. Note also that ${\cal S} = {\cal S}_{\,\overline{1}}$. 

It was stated in~\cite[Observation 2]{goddard-2012} that every graph $G$ and any edge $e$ of it satisfy the inequality 
\begin{align*} 
    \chi_S(G-e) \leq \chi_S(G).
\end{align*}
As the removal of isolated vertices does not change $\chi_S(G)$, this inequality also implies $\chi_S(H) \le \chi_S(G)$ for every subgraph $H$ of $G$. We may also infer that if $\chi_S(G)=k$, then  $G$ contains a subgraph that is $k$-$\chi_S$-critical.

\begin{observation}
\label{obs:edge-removal}
    Let $S \in \cS$ and let $G$ be a graph. 
    \begin{itemize}
        \item[(i)] If $G$ contains no isolated vertex, then $G$ is $\chi_S$-critical if and only if $\chi_S(G-e) < \chi_S(G)$ holds for every edge $e \in E(G)$. 
        \item[(ii)] $K_1$ is the unique $1$-$\chi_S$-critical graph.
    \end{itemize}
    \end{observation}

\section{Families of $S$-packing critical graphs}
\label{sec:families}

In this section we present several families of $S$-packing critical graphs. We first show that graphs of diameter $k$ and girth at least $k+2$ are $\chi_S$-critical for each packing sequence $S \in {\cal S}_k$. This result is then extended to specific generalized lexicographic products. We end the section by classifying $k$-$\chi_S$-critical graphs for almost all $S$ and $k\in \{2,3,4\}$. But first we give two general, simple properties of $S$-packing critical graphs.

\begin{lemma}
\label{lem:connected}
If $S\in {\cal S}$ and $G$ is a $\chi_S$-critical graph, then $G$ is connected.
\end{lemma}
\proof Suppose to the contrary that $G$ is not a connected graph such that $H_1, \dots, H_r$ are the components of $G$, where $r \geq 2$. Since
$
\chi_S(G) = \max_{i \in [r]} \chi_S(H_i),
$
there exists a component $H_j$ such that $\chi_S(H_j) = \chi_S(G)$. Now consider a component $H_k$ for some $k \neq j$. The removal of $H_k$ from $G$ yields a proper subgraph $H$ with
$
\chi_S(H) = \chi_S(H_j) = \chi_S(G),
$
which contradicts the assumption that $G$ is $\chi_S$-critical. Therefore, $G$ must be connected.
\qed

\begin{lemma}
\label{lem:edge-vertex-critical}
If $G$ is a $\chi_S$-critical graph, then $G$ is a $\chi_S$-vertex-critical graph.
\end{lemma}

\proof
If $G \cong K_1$, then it is $\chi_S$-critical  and $\chi_S$-vertex-critical. Otherwise, consider an arbitrary vertex $x \in V(G)$. Since $G-x$ is a proper subgraph of $G$, $\chi_S$-criticality implies $\chi_S(G - x) < \chi_S(G)$ and proves that $G$ is a $\chi_S$-vertex-critical graph. \qed 

A graph $G$ is called {\em diameter $k$-critical} if $\diam(G)=k$ and $\diam(G-e) > \diam(G)$ holds for every $e\in E(G)$ (see~\cite{devillers-2024+, loh-2016, royle-2002}).

\begin{proposition}
\label{prop:diam-critical}
Let $k\ge 1$ and $S \in {\cal S}_k$. Then every diameter $k$-critical graph is $\chi_S$-critical. 
\end{proposition}

\proof
Let $G$ be a diameter $k$-critical graph. Since $s_1 = k$, no two vertices of $G$ can receive the same color, that is, $\chi_S(G) = n(G)$. Let now $e$ be an arbitrary edge of $G$. Since $\diam(G-e) \ge k+1$, there are two vertices $u$ and $v$ with $d_{G-e}(u,v) \ge k+1$. Therefore, in an $S$-coloring of $G-e$, we can color $u$ and $v$ with color $1$, and assign a unique color to every other vertex. Hence $\chi_S(G-e) \le n(G) - 1 < \chi_S(G) = n(G)$ which yields the conclusion. 
\qed

Since every graph $G$ with $\diam(G)=k$ and girth $g(G) \ge k+2$ is diameter $k$-critical, we deduce the following statement from Proposition~\ref{prop:diam-critical}. 
\begin{corollary}
\label{cor:girth-k+2}
Let $k\ge 1$ and $S \in {\cal S}_k$. If $G$ is a graph with $\diam(G) = k$ and $g(G) \ge k+2$, then $G$ is $\chi_S$-critical. 
\end{corollary}


As Proposition~\ref{prop:diam-critical} and Corollary~\ref{cor:girth-k+2} are true for trees, we may infer that every tree is $\chi_S$-critical for infinitely many packing sequences $S$. We also prove the following property for trees.

\begin{proposition}
\label{prop:critical-tree}
   For every tree $T$ and every $S\in {\cal S}$, the tree $T$ is $\chi_S$-critical if and only if it is $\chi_S$-vertex-critical. 
\end{proposition}

\proof
If $T$ is an isolated vertex, the equivalence holds. From now on, we may assume $n(T) \ge 2$.
The first direction of the statement follows from Lemma~\ref{lem:edge-vertex-critical}. To prove the other direction, consider a $\chi_S$-vertex-critical tree $T$ and remove an arbitrary edge $e=u_iu_j$. Let $T_i$ and $T_j$ be the two components of $T-e$ such that $u_i \in V(T_i)$ and $u_j \in V(T_j)$. We know that $\chi_S(T-e)= \max\{\chi_S(T_i), \chi_S(T_j)\}$ and may assume $\chi_S(T-e)= \chi_S(T_i)$. As $T_i$ is also a component in $T-u_j$, the $\chi_S$-vertex-criticality of $T$ implies 
$$
\chi_S(T-e)= \chi_S(T_i) \leq \chi_S(T-u_j) < \chi_S(T).
$$
Since $\chi_S(T-e) < \chi_S(T)$ holds for every edge and $T$ is isolate-free, we may conclude that $T$ is $\chi_S$-critical. This finishes the proof of the equivalence.
\qed

We conclude the section by considering $k$-$\chi_S$-critical graphs, where $k\in \{2,3,4\}$.

\begin{proposition}
\label{prop:2chiS_edge}
If $S\in {\cal S}$, then a graph $G$ is $2$-$\chi_S$-critical if and only if $G \cong K_2$.
\end{proposition}

\proof
It is straightforward that $G \cong K_2$ is $2$-$\chi_S$-critical: its packing chromatic number is $2$, and removing any vertex or edge reduces the packing chromatic number to $1$. Conversely, suppose $G$ is $2$-$\chi_S$-critical. Then, by Lemma~\ref{lem:edge-vertex-critical}, $G$ is also $2$-$\chi_S$-vertex-critical. It is shown in~\cite{holub-2020} that the only graph with this property is $K_2$. Thus, the result follows.
\qed

In the first item of~\cite[Theorem 5.1]{holub-2020} it is claimed that if $S\in {\cal S}_{2,2,2}$, then a graph $G$ is $4$-$\chi_S$-vertex-critical if and only if $G$ is one of the graphs $K_{1,3}$, $C_4$, $Z_1$,  
$K_4-e$, $K_4$, where $Z_1$ denotes the graph obtained by adding a pendant edge to a $C_3$. However, there is one example missing from the proof, which we explain below.

Consider a $4$-$\chi_S$-vertex-critical graph $G$ and let $u\in V(G)$, so that $\chi_S(G-u)\le 3$. In the proof of~\cite[Theorem 5.1]{holub-2020}, when the case $\deg_G(u) = 2$ is considered, it is correctly stated that if $G-u$ is connected, then, for $G-u\cong P_3$ we get $G\cong C_4$ or $G\cong Z_1$, for $G-u\cong C_3$ we get $G\cong K_4-e$. Afterwards, it is stated that in any other case no $\chi_S$-critical graph is obtained. But the vertex $u$ can be adjacent to two degree $1$ vertices of $G-u$, that is, to the end-vertices of a path in which case $G$ is a cycle. Note that $\chi_S(C_n) = 3$ if and only if $n$ is divisible by $3$. Hence, $C_n$ is a $4$-$\chi_S$-vertex-critical graph if and only if $n$ is not divisible by $3$. 

According to the above, the first item of~\cite[Theorem 5.1]{holub-2020} must be supplemented as follows.

\begin{proposition} 
\label{prop:corrected}
If $S \in {\cal S}_{2,2,2}$, then a graph $G$ is $4$-$\chi_S$-critical if and only if 
$$G \in \{K_{1,3}, Z_1, K_4-e, K_4\} \cup \{C_n:\ n\ge 4,\ n\not\equiv 0 \pmod 3\}.$$
\end{proposition} 

With Proposition~\ref{prop:corrected} in hand, we can state the following result. 

\begin{theorem}
\label{thm:many-cases}
Let $S\in \cS$ and let $G$ be a graph. 
\begin{enumerate}
\item[(i)] If $S \in \cS_{1,1}$, then $G$ is $3$-$\chi_S$-critical if and only if $G \in \{C_{2k+1}:k\geq 1\}$.
\item[(ii)] If  $S \in \cS_{1, \overline{2}}$, then $G$ is $3$-$\chi_S$-critical if and only if $G \in \{C_3,P_4\}$
\item[(iii)] If $S \in \cS_{\overline{2}}$, then $G$ is $3$-$\chi_S$-critical if and only if $G \cong P_3$.
\item[(iv)] If $S \in \cS_{2,2,2}$, then $G$ is $4$-$\Spch$-critical if and only if \\ $G \in \{K_{1,3}\} \cup \{C_n:\ n\ge 4,\ n\not\equiv 0 \pmod 3\}$.
\item[(v)] If $S \in \cS_{2,2, \overline{3}}$, then $G$ is $4$-$\Spch$-critical if and only if $G\in \{K_{1,3},C_4, P_6\}$.
\item[(vi)] If $S \in \mathcal{S}_{2, \overline{3}}$, then $G$ is $4$-$\Spch$-critical if and only if $G\in \{K_{1,3},C_4, P_5\}$.
\item[(vii)] If $S \in {\cS}_{\overline{3}}$, then $G$ is $4$-$\Spch$-critical if and only if $G \in \{  K_{1,3}, P_4\}$.
\end{enumerate}
\end{theorem} 

\proof
Let $S\in {\cal S}$. Then it was proved in \cite[Theorem 4.1]{holub-2020} that (i) if  $S \in \cS_{1,1}$, then $G$ is $3$-$\chi_S$-vertex-critical if and only if $G \in \{C_{2k+1}:k\geq 1\}$; (ii) if $S \in \cS_{1, \overline{2}}$, then $G$ is $3$-$\chi_S$-vertex-critical if and only if $G \in \{C_3,C_4,P_4\}$; and (iii) if  $S \in \cS_{\overline{2}}$, then $G$ is $3$-$\chi_S$-vertex-critical if and only if $G \in \{C_3,P_3\}$. By Lemma~\ref{lem:edge-vertex-critical}, it remains to verify which of the above-listed graphs is  $3$-$\chi_S$-critical. Doing it one by one, the first three assertions of the theorem follow.

By Proposition~\ref{prop:corrected} we know the list of $4$-$\Spch$-vertex-critical where $S \in \cS_{2,2,2}$. We next recall that it was further proved in Theorem~\cite[Theorem 5.1]{holub-2020} that (a) if $S \in \cS_{2,2, \overline{3}}$, then $G$ is $4$-$\Spch$-vertex-critical if and only if $G\in \{K_{1,3},C_4,Z_1, K_4-e, K_4, P_6,C_6\}$; (b) if $S \in \cS_{2, \overline{3}}$, then $G$ is $4$-$\Spch$-vertex-critical if and only if  $G\in \{K_{1,3}, C_4, Z_1, K_4-e, K_4,P_5\}$; and (c) if $S \in \cS_{\overline{3}}$, then $G$ is $4$-$\Spch$-vertex-critical if and only if $G \in \{  K_{1,3}, P_4, C_4, Z_1, K_4-e, K_4\}$. Applying Lemma~\ref{lem:edge-vertex-critical} again, we need to verify which of the above-listed graphs are  $4$-$\chi_S$-critical. Carefully checking all of them, the last four assertions of the theorem follow.
\qed

\begin{figure}[ht!]
\begin{center}
\begin{tikzpicture}[scale=1.0,style=thick,x=1cm,y=1cm]
\def\vr{3pt}

\begin{scope}[xshift=0cm, yshift=0cm] 
\coordinate(x1) at (0,0);
\coordinate(x2) at (0,1);
\coordinate(x3) at (1,0.5);
\coordinate(x4) at (2,0.5);
\coordinate(x5) at (3,0.5);
\draw (x1) -- (x2) -- (x3) -- (x1);
\draw (x3) -- (x4) -- (x5);
\foreach \i in {1,...,5}
{
\draw(x\i)[fill=white] circle(\vr);
}
\node at (1.5,-0.5) {$G_1$};
\end{scope}

\begin{scope}[xshift=5cm, yshift=0cm] 
\coordinate(x1) at (0,0);
\coordinate(x2) at (0,1);
\coordinate(x3) at (1,0);
\coordinate(x4) at (1,1);
\coordinate(x5) at (2,0);
\coordinate(x6) at (2,1);
\draw (x5) -- (x3) -- (x1) -- (x2) -- (x4) -- (x6);
\draw (x3) -- (x4);
\foreach \i in {1,...,6}
{
\draw(x\i)[fill=white] circle(\vr);
}
\node at (1,-0.5) {$G_2$};
\end{scope}

\begin{scope}[xshift=9cm, yshift=0cm] 
\coordinate(x1) at (0,0);
\coordinate(x2) at (0,1);
\coordinate(x3) at (1,0);
\coordinate(x4) at (1,1);
\coordinate(x5) at (2,0);
\draw (x5) -- (x3) -- (x1) -- (x2) -- (x4);
\draw (x1) -- (x4) -- (x3);
\foreach \i in {1,...,5}
{
\draw(x\i)[fill=white] circle(\vr);
}
\node at (1,-0.5) {$G_3$};
\end{scope}

\begin{scope}[xshift=1cm, yshift=-2.5cm] 
\coordinate(x1) at (0,0);
\coordinate(x2) at (1,0);
\coordinate(x3) at (2,0);
\coordinate(x4) at (3,0);
\coordinate(x5) at (1.5,0.7);
\coordinate(x6) at (1.5,1.4);
\draw (x1) -- (x2) -- (x3) -- (x4);
\draw (x2) -- (x5) -- (x3);
\draw (x5) -- (x6);
\foreach \i in {1,...,6}
{
\draw(x\i)[fill=white] circle(\vr);
}
\node at (1.5,-0.5) {$G_4$};
\end{scope}

\begin{scope}[xshift=6cm, yshift=-2.5cm] 
\coordinate(x1) at (0,0);
\coordinate(x2) at (1,0);
\coordinate(x3) at (2,0);
\coordinate(x4) at (3,0);
\coordinate(x5) at (4,0);
\coordinate(x6) at (0,1);
\coordinate(x7) at (3,1);
\coordinate(x8) at (4,1);
\draw (x1) -- (x2) -- (x3) -- (x4) -- (x5) -- (x8);
\draw (x1) -- (x6);
\draw (x4) -- (x7);
\foreach \i in {1,...,8}
{
\draw(x\i)[fill=white] circle(\vr);
}
\node at (2.0,-0.5) {$G_5$};
\end{scope}

\begin{scope}[xshift=0cm, yshift=-5cm] 
\coordinate(x1) at (0,0);
\coordinate(x2) at (1,0);
\coordinate(x3) at (2,0);
\coordinate(x4) at (3,0);
\coordinate(x5) at (-0.3,1);
\coordinate(x6) at (0.3,1);
\coordinate(x7) at (2.7,1);
\coordinate(x8) at (3.3,1);
\draw (x5) -- (x1) -- (x2) -- (x3) -- (x4) -- (x8);
\draw (x1) -- (x6);
\draw (x4) -- (x7);
\foreach \i in {1,...,8}
{
\draw(x\i)[fill=white] circle(\vr);
}
\node at (1.5,-0.5) {$G_6$};
\end{scope}

\begin{scope}[xshift=4.5cm, yshift=-5cm] 
\coordinate(x1) at (0,0);
\coordinate(x2) at (1,0);
\coordinate(x3) at (2,0);
\coordinate(x4) at (3,0);
\coordinate(x5) at (0,1);
\coordinate(x6) at (1,1);
\coordinate(x7) at (2,1);
\draw (x5) -- (x1) -- (x2) -- (x3) -- (x4);
\draw (x2) -- (x6);
\draw (x3) -- (x7);
\foreach \i in {1,...,7}
{
\draw(x\i)[fill=white] circle(\vr);
}
\node at (1.5,-0.5) {$G_7$};
\end{scope}

\begin{scope}[xshift=9.0cm, yshift=-5cm] 
\coordinate(x1) at (0,0);
\coordinate(x2) at (1,0);
\coordinate(x3) at (2,0);
\coordinate(x4) at (0,1);
\coordinate(x5) at (1,1);
\coordinate(x6) at (2,1);
\draw (x4) -- (x1) -- (x2) -- (x3) -- (x6);
\draw (x2) -- (x5);
\foreach \i in {1,...,6}
{
\draw(x\i)[fill=white] circle(\vr);
}
\node at (1,-0.5) {$G_8$};
\end{scope}

\end{tikzpicture}
\caption{The graphs $G_1, \ldots, G_8$}
\label{fig:sporadic-graphs}
\end{center}
\end{figure}
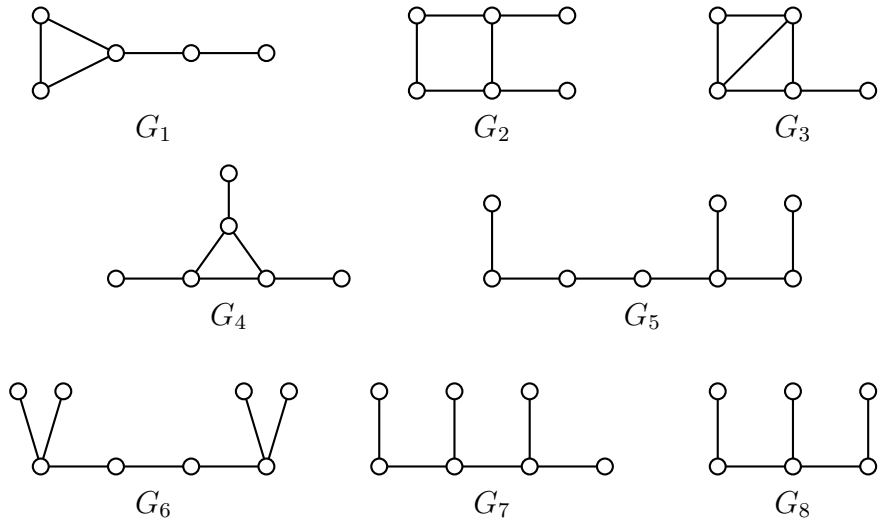

To state the following theorem, we first need to introduce two families of graphs and some specific graphs. If $S \in \cS$, then let 
$$\mathcal{C}_{s_4} = \{C_n, n \ge 5:\  (n\equiv 1,2 \pmod{4}) ~\text{or}~ (n\equiv3 \pmod{4}~\text{and}~s_4<\lfloor n/2 \rfloor) \}\,.$$
So $\mathcal{C}_{s_4}$ is a subclass of cycles that depends on the fourth term of $S$. Next, let $X_{2k}$, $k\ge 3$, be the graph obtained from $P_{2k}$ by attaching a pendant vertex to each of the support vertices of $P_{2k}$. Finally, we need the graphs $G_i$, $i\in \{1,\dots, 8\}$, which are shown in Fig.~\ref{fig:sporadic-graphs}. Note that $G_6 \cong X_6$.

\begin{theorem}
\label{thm:4chiS_smallcases}
If $G$ is a graph, then the following assertions hold.
\begin{enumerate}

\item[(i)] If $S\in \mathcal{S}_{1,3,3}$, then $G$ is $4$-$\Spch$-critical if and only if \[G \in \{K_4, G_1, G_2\} \cup \mathcal{C}_{s_4}  \cup \{X_{2k}: k\geq 3\}.\]
\item[(ii)] If $S\in \mathcal{S}_{1,3,\overline{4}}$, then $G$ is $4$-$\Spch$-critical if and only if \[G \in \{K_4,C_5,C_6,P_8\} \cup \{G_i:\ i\in \{1,\dots, 7\}\}.\]
\item[(iii)] If $S\in \mathcal{S}_{1,\overline{4}}$, then $G$ is $4$-$\Spch$-critical if and only if \[G \in \{K_4, C_5, P_6, G_8\}.\]
\end{enumerate}
\end{theorem}
\proof 
Let $S\in {\cal S}$. In each of the three cases listed above, the set of all $4$-$\chi_S$-vertex-critical graphs has been completely classified in~\cite{klavzar-2023}. 
This classification includes sporadic examples, as well as several graph families.  
By Lemma~\ref{lem:edge-vertex-critical}, any $4$-$\chi_S$-critical graph must also be $4$-$\chi_S$-vertex-critical.  
Therefore, it suffices to determine which of the graphs identified in~\cite{klavzar-2023} are also $\chi_S$-critical.

When $S \in \mathcal{S}_{1,3,3}$, the classification from~\cite{klavzar-2023} includes nine individual graphs and four graph families, two of which are infinite.  
Only those listed in item~(i) satisfy the $\chi_S$-criticality condition.

For $S \in \mathcal{S}_{1,3,\overline{4}}$, the classification from~\cite{klavzar-2023} consists of thirteen graphs and two graph families.  
Only the graphs listed in item~(ii) are $\chi_S$-critical.

For $S \in \mathcal{S}_{1,\overline{4}}$, the classification from~\cite{klavzar-2023} includes twelve graphs and two graph families.  
Among them, only the graphs in item~(iii) were found to be $\chi_S$-critical through direct verification. This completes the proof.
\qed

\section{On $S$-packing critical cycles}
\label{sec:cycles}

Recall that $(k,k,\ldots)$-packing colorings are known as $k$-distance colorings and that the $(k,k,\ldots)$-packing chromatic number is denoted by $\chi_k$. For cycles, it is proved in~\cite{naranjan-2019} that  $\chi_k(C_n) = k+1 + \lceil \frac{r}{\ell} \rceil$, where $n=\ell(k+1) + r$, $0\le r\le k$. (The special case when $r=0$ was earlier observed in~\cite{kramer-1969}.) Since $n\ge k+1$ implies  $\chi_k(P_n) = k+1$, we can conclude that $C_n$, where $n\ge k+1$, is $(k,k,,\ldots)$-packing critical if and only if $n \not\equiv 0 \pmod {k+1}$.

\begin{proposition}
\label{prop:cycles-small-n}
If $S\in {\cal S}$ and $n\ge 3$, then the following hold. 
\begin{enumerate}
\item[(i)] If $n \le s_1 + 1$, then $C_n$ is not $S$-packing critical.
\item[(ii)] If $s_1+2 \le n \le 2s_1 + 1$, then $C_n$ is $S$-packing critical.
\end{enumerate}
\end{proposition}

\proof
(i) Since $n \le s_1 + 1$, we have $\chi_S(C_n) = n$. Moreover, $C_n-e = P_n$ and $\diam(P_n) = n-1 \le s_1$, hence $\chi_S(C_n-e) = n$. So $C_n$ is not $S$-packing critical.

(ii) Let $s_1+2 \le n \le 2s_1 + 1$. Since $\diam(C_n) = \lfloor \frac{n}{2}\rfloor \le \lfloor \frac{2s_1+1}{2}\rfloor = s_1$, we have $\chi_S(C_n) = n$. On the other hand, $\diam(C_n-e) = \diam(P_n) = n-1 \ge (s_1+2) - 1 = s_1+1$, so we can color the two leaves of $P_n$ with the same color and henceforth, $\chi_S(C_n-e) \le n-1$. So $C_n$ is $S$-packing critical in this case.
\qed

\begin{theorem}
\label{prop:cycles-1-big}
If $S \in {\cal S}_{1}$ and $n\ge 3$, then the following hold. 
\begin{enumerate}
\item[(i)] If $n \le s_2 + 2$, then $C_n$ is  $S$-packing critical if and only if $n$ is odd.
\item[(ii)] If $s_2+3 \le n \le 2s_2 + 1$, then $C_n$ is $S$-packing critical.
\end{enumerate}
\end{theorem}

\proof
(i) If $n \leq s_2+1$, then $\diam(C_n) < \diam(P_n) \le s_2$. Therefore, no color $t$ with $t \ge 2$ can be repeated in an $S$-packing coloring of $P_n$, and we get an optimal coloring by assigning color $1$ to as many vertices as possible. Then $\chi_S(P_n)= n- \alpha(P_n)+1= n-\lceil\frac{n}{2}\rceil+1$. For cycles, the situation is similar and we obtain $\chi_S(C_n)= n- \alpha(C_n)+1= n-\lfloor\frac{n}{2}\rfloor+1$. Since $\chi_S(C_n -e) = \chi_S(P_n)$ for every $e \in E(C_n)$, the above formulas show that $\chi_S(C_n -e)  < \chi_S(C_n)$ if and only if $n$ is odd.

Now let $n=s_{2}+2$. We observe that $\diam(C_n) =\lfloor \frac{n}{2} \rfloor \le s_2$ always holds and then no color $t$ with $t \ge 2$ can be repeated in an $S$-packing coloring of $C_n$. Thus, 
$\chi_{S}(C_n) = n-\alpha(C_n) +1 = n- \left\lfloor \frac{n}{2} \right\rfloor +1 = \left\lceil \frac{n}{2} \right\rceil +1.
$ After deleting an edge from $C_n$, we obtain $P_n$, whose diameter is $n-1 =s_2 +1$. We may have two types of $S$-coloring $c$ of the $n$-path $v_1 \dots v_n$. The first possibility is that no color $t$ with $t \ge 2$ is repeated, and then $c$ uses at least $n-\alpha(P_n) +1= \lfloor \frac{n}{2} \rfloor +1 $ colors. The second possibility is to assign color $2$ to $v_1$ and $v_n$, and use color $1$ on an independent set in $v_2 \dots v_{n-1}$. This needs at least $n-2- \alpha(P_{n-2}) +2= n - \lceil \frac{n-2}{2} \rceil = \lfloor \frac{n}{2} \rfloor +1$ colors. We may therefore infer that $\chi_S(C_n) = \lfloor \frac{n}{2} \rfloor +1$. Comparing $\chi_S(C_n)$ and $\chi_S(C_n-e)=\chi_S(P_n)$, we conclude that $C_n$ is $\chi_S$-critical if and only if $n$ is odd, as stated.

(ii) Suppose first that $n$ is even. Since $\diam(C_n) = \frac{n}{2} \le s_2$, no color $t$ with $t \ge 2$ can be repeated in an $S$-packing coloring of $C_n$. Thus,  $\chi_S(C_n)= n- \alpha(C_n)+1= \frac{n}{2} +1.$ For the $n$-path $v_1v_2\dots v_n$, consider the coloring $c$ that assigns color $1$ to the vertices $v_1, v_3, \dots , v_{n-1}$, color $2$ to $v_2$ and $v_n$, while the remaining vertices are colored pairwise differently with colors $3, \dots, \frac{n}{2}$. As $d_{P_n}(v_2,v_n) =n-2 \ge s_2+1$, $c$ is an $S$-packing-coloring. Then, we conclude that $\chi_S(P_n) \leq \frac{n}{2} < \chi_S(C_n)$, proving the $ \chi_S$-criticality of $C_n$.

If $n$ is odd, $\diam(C_n) = \frac{n-1}{2} \leq  s_2$ implies that no color different from $1$ can be repeated in an $S$-packing coloring. We infer again that $\chi_S(C_n)= n- \alpha(C_n)+1= \frac{n+1}{2} +1.$ The path $P_n$ can be colored such that only color $1$ is repeated. Hence,
$$ \chi_S(P_n) \leq  n- \alpha(P_n)+1= \frac{n-1}{2}+1 <  \chi_S(C_n)
$$
that proves the $ \chi_S$-criticality of $C_n$.
\qed

\begin{theorem}
\label{thm:cycles-1-small}
If $n\ge 3$, then the following hold. 
\begin{enumerate}
\item[(i)] If $S \in \mathcal{S}_{1,1}$, then $C_n$ is  $S$-packing critical if and only if $n$ is odd.
\item[(ii)] If $S \in \mathcal{S}_{1,2,2}$, then $C_n$ is  $S$-packing critical if and only if it is $C_3$ or $C_5$.
\item[(iii)] If $S \in \mathcal{S}_{1,\overline{2},3}$, then $C_n$ is  $S$-packing critical if and only if $n \not\equiv 0 \pmod 4$.
\end{enumerate}
\end{theorem}

\proof
Throughout the proof, let $v_1,\dots, v_n$ be consecutive vertices of $C_n$. 

(i) Let $S\in \mathcal{S}_{1,1}$. An even cycle $C_n$ can be colored alternately with colors $1$ and $2$. Hence, $\chi_S(C_n)=2$. If $n$ is odd, a $2$-packing-coloring is not possible, but three colors are clearly enough. On the other hand, $\chi_S(P_n) =2$ for every $n \ge 2$. It follows that $\chi_S(C_n-e) < \chi_S(C_n)$ holds if and only if $n$ is odd.

(ii) Let $S\in \mathcal{S}_{1,2,2}$. Consider a path $P_n$, for $n \ge 4$, and an $S$-packing-coloring $c$ of $P_n$. On every four consecutive vertices of the path, the coloring uses at least three colors. Let $(123)^\ast$ denote the sequence of colors in which $123$ repeats an arbitrary number of times. Using the color pattern $(123)^\ast$, starting with the first vertex of the path, and where from the last block $123$ the required number of elements is used (possibly zero), we obtain an $S$-packing coloring. Consequently, $\chi_S(P_n)=3$ if $n \ge 4$. As follows, $\chi_S(C_n) \ge 3$ holds for every $n \ge 4$. 

If $n \ge 6$,
we consider the following colorings of $C_n$. The referred patterns start from vertex $v_1$, and after a specified initial sequence, the coloring repeats pattern $123$ so that the color of $v_n$ will be $3$. If $n\equiv 0 \pmod 3$, we color $C_n$ with $(123)^\ast$. If $n\equiv 1 \pmod 3$, we color $C_n$ as $1213\,(123)^\ast$. If $n\equiv 2 \pmod 3$, then $n \ge 8$, and we can color $C_n$ as $1213\,1213\,(123)^\ast$. It shows $\chi_S(C_n)\leq 3$ and in turn, $\chi_S(C_n)= 3$, for every $n \ge 6$. We conclude that in this case there is no $S$-packing critical cycle on more than $5$ vertices.

For the small cases, we observe 
$\chi_S(P_3) = 2 < \chi_S(C_3) =3$;
$\chi_S(P_4) = 3 = \chi_S(C_4)$; and 
$\chi_S(P_5) = 3 < \chi_S(C_5) =4$. Now, we may conclude that $C_3$ and $C_5$ are the only $S$-packing critical cycles if $S \in \cS_{1,2,2}$.

(iii) Under the conditions $n \ge 4$ and $S\in \mathcal{S}_{1,\overline{2},3}$, any $S$-packing-coloring of $P_n$ or $C_n$ requires at least $3$ colors. A path $P_n$ with $n \ge 4$ can be colored by $(1213)^\ast$ no matter whether $n \equiv 0 \pmod 4$ is valid or not. Naturally, if $n \not\equiv 0 \pmod 4$, then from the last block $1213$ the required number of elements is used. Therefore, $\chi_S(P_n)=3$ when $n \ge 4$. If $n \equiv 0 \pmod 4$, we can take the same type of coloring for $C_n$ and get $\chi_S(C_n)=3$. It also shows that no $n$-cycle with $n \equiv 0 \pmod 4$ is $\chi_S$-critical.

Suppose now that $c$ is an $S$-packing coloring of $C_n$, for $n \ge 5$ that uses only colors $1$, $2$, $3$. We claim that there are no two neighbors colored with $2$ and $3$. Assume, without loss of generality, that $c(v_i)=3$ and $c(v_{i+1})=2$. Then, $v_{i+2}$ cannot get a color different from $1$. But then, as $s_2\ge 2$ and $s_3 =3$, neither of colors $1$, $2$, and $3$ can be assigned to  $v_{i+3}$. This contradiction proves that every second vertex of the cycle is colored with $1$. As neither of the patterns $1212$ and $1313$ may occur in the coloring, we obtain that the pattern $1213$ must be repeated along the cycle. If $n \not\equiv 0 \pmod 4$, it is impossible to have $3$ colors and we conclude $\chi_S(C_n) \ge 4$ for these cases.  Therefore, $C_n$ is $\chi_S$-critical for every $n \ge 5$ if $n \not\equiv 0 \pmod 4$. Observing also that $\chi_S(P_3)=2 < \chi_S(C_3)=3$ we obtain that $C_3$ is $\chi_S$-critical. This completes the proof for (iii).
\qed

\section{Impact of edge removal on $\chi_S$}
\label{sec:edge-removal}

In view of Observation~\ref{obs:edge-removal}~(i), the question naturally arises as to what extent removing an edge of $G$ can affect $\chi_S(G)$. Before we answer this question, recall the following well-known sets (see~\cite{ali-2021, jerebic-2008, miklavic-2018}) which are defined for an arbitrary edge $e=uv$ of a graph $G$:
\begin{align*}
W_{uv}^{G} & = \{w\in V(G):\ d_G(u,w) < d_G(v,w) \}, \\
W_{vu}^{G} & = \{w\in V(G):\ d_G(v,w) < d_G(u,w) \}, \\
_vW_u^{G} & = \{w\in V(G):\ d_G(u,w) = d_G(v,w) \}\,.
\end{align*}
Clearly, $V(G)= W_{uv} \cup W_{vu} \cup\, _vW_u$. We will use the next lemma throughout the rest of the section mostly without explicitly mentioning it. 

\begin{lemma}
\label{lem:Wuv-edge-removed}
If $e=uv\in E(G)$, then $W_{uv}^{G} = W_{uv}^{G-e}$ and $W_{vu}^{G} = W_{vu}^{G-e}$.
\end{lemma} 

\proof
Assume first that $w\in W_{uv}^{G}$. Then $e=uv$ does not lie on any shortest $w,u$-path, thus we have 
$$d_{G-e}(w,u) = d_{G}(w,u) < d_{G}(w,v) \le d_{G-e}(w,v)\,,$$
hence $w\in W_{uv}^{G-e}$. 

Assume second that $w\in W_{uv}^{G-e}$, that is, $d_{G-e}(w,u) < d_{G-e}(w,v)$. Then no matter whether there exists a shortest $w,v$-path in $G$ which passes $e$, we have 
$$d_{G}(w,v) \ge d_{G}(w,u) + 1\,,$$
that is, $w\in W_{uv}^{G}$. We can conclude that $W_{uv}^{G} = W_{uv}^{G-e}$. The argument for the equality $W_{vu}^{G} = W_{vu}^{G-e}$ is parallel.
\qed

In the proof of the next result, we use some ideas similar to those in the proof of~\cite[Theorem 1]{bresar-2022}.

\begin{theorem}
\label{thm:edge-removed-general}
Let $S$ be a packing sequence and let $e=uv$ be an edge of a graph $G$. Then the following statements hold.
\begin{enumerate}
\item[(i)] $\displaystyle{ \chi_S(G-e) \ge \frac{\chi_S(G)}{2}}$. Moreover, there are infinitely many sharp examples for every packing sequence $S \in  \cS_{1, \overline{3}} \cup \cS_{2, \overline{5}} \cup \cS_{\, \overline{3}}$.
\item[(ii)] If $G$ contains a component on at least three vertices and $S \in \cS_{1,1} \cup \cS_{1,2}$, then $\displaystyle{\chi_S(G-e) \ge \frac{\chi_S(G)+1}{2} }$ holds. 
 Moreover, there are infinitely many sharp examples for every  $S \in \cS_{1,1} \cup \cS_{1,2}$.
 \item[(iii)] If $S \in \cS_{2,2,2}$ and $\chi_S(G-e) \ge 3$, then $\displaystyle{\chi_S(G-e) \ge \frac{\chi_S(G)+1}{2} }$ holds. 
\end{enumerate}
\end{theorem}

\proof
(i)  Let $c'\colon V(G) \rightarrow [\chi_S(G-e)]$ be a $\chi_S$-packing-coloring of $G-e$. For a color $t \in [\chi_S(G-e)]$, we say that a vertex pair $\{x,y\}$ is \emph{$t$-problematic} if $c'(x)=c'(y)=t$ but $d_G(x,y) \leq s_t$. Since $c'$ is an $S$-packing coloring of $G-e$, we have $d_{G-e} (x,y) \ge s_t +1$. Then $d_{G-e} (x,y) > d_G(x,y)$ and therefore, in $G$, every shortest $(x,y)$-path goes through the edge $e$. It also follows that, for every problematic pair $\{x,y\}$, one vertex is in $W_{uv}^G$ and the other is in $W_{vu}^G$. Note that $W_{uv}^{G} = W_{uv}^{G-e}$ and $W_{vu}^G = W_{vu}^{G-e}$ hold by Lemma~\ref{lem:Wuv-edge-removed}. 
\medskip

We say that a vertex $z$ \emph{covers} a problematic pair $\{x,y\}$ if $z=x$ or $z=y$ and state the following key property of problematic pairs.

\medskip\noindent
{\bf Claim}: For every $t \in [\chi_S(G-e)]$, either there is no $t$-problematic pair or there exists a vertex that covers all $t$-problematic pairs. 
\medskip

\noindent \textit{Proof.} Consider the bipartite graph $F_t$ with partite classes $W_{uv}^G$, $W_{vu}^G$, where $xy$ is an edge if $\{x,y\}$ is a $t$-problematic pair in $G$. Suppose for a contradiction that the claim is not true, that is, $E(F_t) \neq \emptyset$ and that one vertex cannot cover all edges of $F_t$. K\"{o}nig's theorem~\cite{konig} implies that the matching number of $F_t$ is at least $2$. So we may suppose that $\{x_1,y_1\}$ and $\{x_2,y_2\}$ are two vertex-disjoint $t$-problematic pairs in $G$. 

Without loss of generality, let $x_i \in W_{uv}^G$ and $y_i \in W_{vu}^G$ for $i \in \{1, 2\}$. Let us set $d_G(x_i,u)=a_i$ and $d_G(y_i,v)=b_i$ for $i \in \{1, 2\}$. Note that these distances remain the same in $G-e$.  Consider first $x_1$ and $x_2$. As both vertices belong to $W_{uv}$, we have  $d_G(x_1,x_2)= d_{G-e}(x_1,x_2) $. Since $c'$ is a $\chi_S$-packing-coloring of $G-e$, it holds that $d_{G-e}(x_1,x_2) \ge s_t +1$. Further, the length  $a_1+a_2$ of the $(x_1, x_2)$-path through $u$ gives an upper bound on the distance between $x_1$ and $x_2$. We obtain
\begin{equation} \label{eq:1}
    a_1 +a_2 \geq d_G(x_1,x_2) \geq s_t+1 . 
        \end{equation}
A similar reasoning gives 
\begin{equation} \label{eq:2}
    b_1 +b_2 \geq d_G(y_1,y_2) \geq s_t+1.  
        \end{equation}
By our assumption, both $\{x_1, y_1\}$ and $\{x_2, y_2\}$ are $t$-problematic pairs and so
\begin{equation} \label{eq:3}
    a_1 +1 +b_1  = d_G(x_1,y_1)\leq  s_t 
        \end{equation}
        and 
      \begin{equation} \label{eq:4}
    a_2 +1 +b_2  = d_G(x_2,y_2)\leq  s_t. 
        \end{equation}  
Inequalities~\eqref{eq:1}-\eqref{eq:4} imply 
\begin{equation*} 
    2 s_t +2 \leq a_1 +a_2 +b_1 +b_2   \leq  2s_t-2.
        \end{equation*}
This contradiction finishes the proof of the claim. \smallqed
\medskip

By the claim, for every color $t$ with a $t$-problematic pair, we can specify a vertex $z_t$ that covers all $t$-problematic pairs. If we remove $z_t$ from the corresponding color class, then no $t$-problematic pair remains, and hence, any two remaining vertices have a distance of at least $s_t+1$ in $G$. Let $Z$ contain all specified vertices $z_t$. Then $|Z| \leq \chi_S(G-e)$. Define now a new coloring $c$ which keeps the color $c'(x)$ if $x \notin Z$ and assigns a unique color to every vertex $x \in Z$ from $\{\chi_S(G-e) +1, \dots , \chi_S(G-e) + |Z|\} $.  

It is clear that $c$ uses at most $2 \chi_S(G-e)$ colors. We now prove that $c$ is an $ S$-packing coloring of $G$. Suppose that $c(x)=c(y)=p$, where $x\ne y$. Since every color $q$ with $q > \chi_S(G-e)$ is assigned to only one vertex, we infer that $p \in [\chi_S(G-e)]$. As all $p$-problematic pairs were destroyed by recoloring one vertex from the pair, $\{x,y\}$ is not a problematic pair and hence, $d_G(x,y) \geq s_p +1$. Thus, $c$ is an $S$-packing coloring of $G$, which implies $\chi_S(G) \leq 2 \chi_S(G-e)$ as stated.
\medskip

We now prove the sharpness of the inequality. If $S \in \cS_{1, \overline{3}}$, let $G$ be constructed by taking two copies of the star $K_{1,k}$ with $k \ge 3$ and connecting them by an edge $e$ between two leaves. It is clear that $\chi_S(G-e)= \chi_S(K_{1,k})=2$. We show that $\chi_S(G)= 4$. In $G$, the path $P$ between the centers of the stars is an isometric subgraph of diameter $3$. Hence, either all four vertices of $P$ get different colors, or color $1$ is assigned to two vertices. In the latter case, at least one center receives color $1$, and then the $k$ neighbors get pairwise different colors. In either case, the number of colors is at least $4$. On the other hand, a $4$-packing-coloring can be obtained by assigning color $1$ to all leaves and one vertex of degree $2$. Thus, $\chi_S(G)= 4$ and $G$  is a sharp example for the bound in (i). 

If $S \in \cS_{2, \overline{5}}$, let $G_k$, $k \ge 2$, be the graph obtained from the disjoint union of $K_k$ and $K_{k+1}$ by adding a path of length $3$ between a vertex of $K_k$ and a vertex of $K_{k+1}$. Let $e$ be the edge of this path attached to $K_{k+1}$. As $\diam(G_k)=5$, no color except $1$ can be repeated in an $S$-packing coloring of $G_k$ and it is easy to check that $\chi_S(G_k)=2k+2=2 \chi_S(G_k-e)$.

Assume now that $S \in S_{\, \overline{3}}$
 and consider the following example. Let $H$ be a graph with a universal vertex and let $G$ be the graph obtained from the disjoint union of two copies of $H$ by adding an edge $e$ between a universal vertex of one copy of $H$ and a universal vertex of the other copy of $H$. Then $\diam(H) = 3$ which implies that $\chi_S(G) = 2n(H)$. On the other hand,  $\chi_S(G-e) = n(H)$. This demonstrates the sharpness of (i) for every $S \in S_{\, \overline{3}}$.

\medskip

(ii) Let $S \in \cS_{1,1} \cup \cS_{1,2}$. If the largest component of $G$ contains at least three vertices, $\chi_S(G-e) \ge 2$ holds for every $e \in E(G)$.
 We prove that there is a color $t \in \{1,2\}$ without a $t$-problematic pair in $G$. Assume that $\{x,y\}$ is a $1$-problematic pair. Then $d_G(x,y) \leq s_1=1$ and all shortest $(x,y)$-paths contain $e=uv$. It implies $\{x,y\}=\{u,v\}$ and $c'(u)=c'(v)=1$. Consequently, for every two vertices $x'$ and $y'$ with $c'(x')=c'(y')=2$, either $d_G(x',y') = d_{G-e}(x',y') \ge s_2+1$ or, in $G$, every shortest $(x',y')$-path contains $e$ and $d_G(x',y') \ge 3 \ge s_2+1$. It follows that one of the colors $1$ and $2$ has no problematic pair, and then, the proof of part (i) can be improved by claiming $|Z| \leq \chi_S(G-e)-1 $. We conclude $\chi_S(G) \leq 2\chi_S(G-e)-1$ as stated.

 For a packing sequence $S \in \cS_{1,1}$, we take the odd cycles which are $3$-$\chi_S$-critical graphs according to Theorem~\ref{thm:many-cases}~(i). Thus, $\chi_S(C_{2k+1})=3$ and $\chi_S(C_{2k+1}-e)=2$, and the odd cycles are sharp examples for the inequality in (ii).

 When $S \in \cS_{1,2}$, we consider two vertex-disjoint stars $K_{1,k}$, for $k \ge 3$, and add an edge $e$ between the centers to obtain the graph $G$. It is easy to check that $\chi_S(G-e)= 2$ and $\chi_S(G)=3$. It provides then a sharp example for (ii). Remark that $C_3$ and $P_4$ are also sharp examples for $S \in \cS_{1,2}$, according to Theorem~\ref{thm:many-cases}~(ii).

(iii) Assume that $S \in \cS_{2,2,2}$ and $\chi_S(G-e) \ge 3$. We prove that for at least one color $t \in \{1,2,3\}$, $G$ contains no $t$-problematic pair. Let us choose $t$ from $\{1,2,3\}$ such that $t\neq c'(u)$ and $t \neq c'(v).$  Then, for every two vertices $x$ and $y$ with $c'(x)=c'(y)=t$, all shortest $(x,y)$-paths contain $e = uv$ and the distance $d_G(x,y)$ is at least $3=s_t+1$. Therefore, we have $|Z| \leq \chi_S(G-e)-1 $ again and may conclude $\chi_S(G) \leq 2\chi_S(G-e)-1$.
\qed

We note that the inequalities in Theorem~\ref{thm:edge-removed-general} (i) and (ii) remain valid if the packing sequence $S$ is finite and we suppose that $G$ is $S$-packing colorable. Indeed, if $2 \chi_S(G-e) \le |S|$, the proof given above remains valid. If 
$2 \chi_S(G-e) > |S|$, then the $ S$-packing colorability of $G$ immediately implies $\chi_S(G) \leq |S| \leq  2\chi_S(G-e)-1$ and the two inequalities follow. 

Setting $S = (1,2,3,\ldots)$ in Theorem~\ref{thm:edge-removed-general}~(ii), we get the following:

\begin{corollary} {\rm \cite[Theorem 1]{bresar-2022}}
\label{cor-bostjan-jasmina}
If $e\in E(G)$, then $\chi_\rho(G-e) \ge \frac{\chi_\rho(G)+1}{2}$.
\end{corollary}

To see that the bound in Theorem~\ref{thm:edge-removed-general}~(i) is asymptotically sharp also when $e$ is not a cut-edge, consider the following example for the constant packing sequence $S=(3,3, \dots)$. Let $H$ be a graph with two universal vertices $x$ and $y$, and let $H'$ be an isomorphic copy of $H$ with respective universal vertices $x'$ and $y'$. Let $G$ be the graph obtained from the disjoint union of $H$ and $H'$ by adding the edge $e = xx'$, and by connecting $y$ and $y'$ with a path of length $3$. 

Note that $n(G) = 2n(H) + 2$ and that $\diam(G) = 3$. Therefore, $\chi_S(G) = 2n(H) + 2$. Consider now $G-e$. Then we can assign color $1$ to $x$ and $y'$, color $2$ to $y$ and $x'$, whilst assigning each color from $\{3, \dots, n(H)\}$ to the remaining pairs of vertices respectively, one from each of $H$ and $H'$. Two further colors, $n(H)+1$ and $n(H)+2$ are used to color the two vertices outside $V(H) \cup V(H')$. In this way, we infer that $\chi_S(G-e)=  n(H) + 2$. So $\lim_{n\to \infty} \frac{\chi_S(G-e)}{\chi_S(G)} = \frac{1}{2}$. 

\medskip
If the removed edge is a cut-edge, we can slightly improve Theorem~\ref{thm:edge-removed-general}.
\begin{proposition} \label{prop:cut-edge}
Let $S\in {\cal S}$ and $s_2 \leq 2$. If $e$ is a cut-edge in a graph $G$ and $\chi_S(G-e) \ge 2$, then $\displaystyle{\chi_S(G-e) \ge \frac{\chi_S(G)+1}{2} }$.
\end{proposition}

\proof Theorem~\ref{thm:edge-removed-general}~(ii) establishes the lower bound if $s_1=1$ and $s_2 \le 2$. Hence, it suffices to prove the lower bound for $s_1=s_2=2$.
Let $e=uv$ be a cut-edge in $G$, and $G_1$, $G_2$ be the two components in $G-e$. We may suppose that $u \in V(G_1)$ and $v \in V(G_2)$. We use the notations from the proof of Theorem~\ref{thm:edge-removed-general}. Assume first that some color $t \in \{1,2\}$ is not in $\{c'(u), c'(v)\}$ and $c'(x) = c'(y)=t$. If $x$ and $y$ belong to the same component $G_i$, then $d_G(x,y)= d_{G-e}(x,y) \geq s_t+1$ as $c'$ is an $S$-packing coloring in $G-e$. If $x \in V(G_1)$ and $y \in V(G_2)$, then the distance $d_G(x,y) \ge 3 =s_t+1$. We conclude that there is no $t$-problematic pair in $G$ and $\chi_S(G) \le 2\chi_S(G-e)-1 $ holds for this case.

If both colors $1$ and $2$ are used on vertices $u$, $v$ by $c'$, we define a coloring $c''$ of $G-e$ by switching colors $1$ and $2$ in $G_2$. Since $s_1=s_2$, coloring $c''$ remains an $S$-packing coloring. Moreover, as $c''(u)=c''(v)$ holds,  $\chi_S(G) \le 2\chi_S(G-e)-1 $ follows by the same reasoning as above.
\qed

In the sharp examples with $\chi_S(G) = 2 \chi_S (G-e)$ given in the proof of Theorem~\ref{thm:edge-removed-general}~(i), the edge $e$ is always a cut-edge. Therefore, the inequality in Proposition~\ref{prop:cut-edge} does not hold for all graphs when  $S \in  \cS_{1, \overline{3}} \cup \cS_{2, \overline{5}} \cup \cS_{\, \overline{3}}$.

\section{Concluding remarks} \label{sec:concluding}

\begin{itemize}
\item In Theorems~\ref{thm:many-cases} and~\ref{thm:4chiS_smallcases} we have characterized $4$-$\Spch$-critical graphs for most of the packing sequences $S$. The missing cases which remain to be considered are $S\in {\cal S}_{1,1} \cup {\cal S}_{1,2}$. In fact, these are also the missing cases of $4$-$\Spch$-vertex-critical graphs. 

\item In Theorem~\ref{thm:cycles-1-small} we have characterized cycles which are $S$-packing critical for $S \in \mathcal{S}_{1,1}$, $S \in \mathcal{S}_{1,2,2}$, and $S \in \mathcal{S}_{1,\overline{2},3}$. The remaining cases are still to be explored.

\item 
In Theorem~\ref{thm:edge-removed-general} we have demonstrated that there are infinitely many sharp examples for the inequality $\chi_S(G-e) \ge \frac{\chi_S(G)}{2}$, for each $S \in  \cS_{1, \overline{3}} \cup \cS_{2, \overline{5}} \cup \cS_{\, \overline{3}}$, where $e$ is a cut-edge. We next provide another sporadic example for the sharpness when $S \in \cS_{2,3,\overline{11}}$. For this purpose, consider $P_{14}$ and its middle edge $e$. Using case analysis, it can be checked that $\chi_S(P_{14})= 8$. On the other hand, $P_{14}-e$ contains two components both of which are isomorphic to $P_7$ and we obtain $\chi_S(P_{14}-e)= 4= \frac{\chi_S(P_{14})}{2}$.
Proposition~\ref{prop:cut-edge} shows that if $s_2 \leq 2$ and $G$ contains a component with at least two edges, then the stronger inequality $\chi_S(G-e) \ge \frac{\chi_S(G)+1}{2}$ holds for every cut-edge $e$ of $G$. The remaining cases are packing sequences with
\begin{itemize}
    \item[$\circ$] $s_1=2$, $s_2=3$, and $3 \le s_3 \le 10$;
    \item[$\circ$] $s_1 =2$, $s_2 =4$.   
\end{itemize} 
For these cases, it remains an open question whether $\chi_S(G-e) \ge \frac{\chi_S(G)+1}{2}$ holds whenever $e$ is a cut-edge of $G$. 

\item In the above example when $S \in \cS_{2,3,\overline{11}}$, we have stated that  $\chi_S(P_{14})= 8$. Establishing this result is not completely straightforward. In general, it would be of interest to determine $\chi_S(P_{n})$ for any $S\in {\cal S}$ and any $n$.

\end{itemize}

\section*{Acknowledgements}

This work was supported by the Slovenian Research and Innovation Agency ARIS (research core funding P1-0297 and projects N1-0285, N1-0355), by the Scientific and Technological Research Council of Türkiye (TÜBİTAK) under grant no.\ 124F114, and by TÜBITAK 2221 Fellowships for Visiting Scientists and Scientists on Sabbatical Leave program.

\baselineskip13pt

\end{document}